\documentclass[conference,11pt]{IEEEtran}
\usepackage{cite}
\usepackage{graphicx}
\usepackage{epsfig}
\usepackage{balance}
\usepackage{array}
\begin{document}

\title{A Realizable Modified Tent Map for True Random Number Generation}
\author{Hamid Nejati, Ahmad Beirami, and Yehia Massoud\\Electrical and Computer Engineering Department, Rice~University,~Houston~TX~77005\\(massoud@rice.edu)}
\maketitle

\begin{abstract}
Tent map is a discrete-time piecewise-affine I/O characteristic curve, which is used for chaos-based applications, such as true random number generation. However, tent map suffers from the inability to maintain the output state confined to the input range under noise and process variations. In this paper, we propose a modified tent map, which is interchangeable with the tent map for practical applications. In the proposed modified tent map, the confinement problem is solved while maintaining the functionality of the tent map. We also demonstrate simulation results for the circuit implementation of the presented modified tent map for true random number generation.
\end{abstract}

\section{Introduction}
\label{sec:intro}

Certain nonlinear dynamical systems can demonstrate non-periodic, long-term non-predictive behaviors, known as chaos. Chaotic behavior results from the high sensitivity of the system to the initial state, which is never exactly known in practice. Any small perturbation can grow exponentially with time within the system leading to a non-predictive behavior of chaotic systems~\cite{stavroulakis06,chaos_dynamic}.
Chaotic waveforms have been extensively used in various research areas such as the modeling of the behavior of human organs~\cite{bio_chaos,tsuda01,chaos_circuit} as well as the modulation of signals and chaotic encryption of telecommunication data~\cite{stavroulakis06}.
Truly Random Number Generators (TRNG) can be utilized in cryptographic systems.
A TRNG is a number generator that is capable of producing uncorrelated and unbiased binary digits through a nondeterministic and irreproducible process~\cite{callegari05a}. A TRNG requires a high entropy source, which can be provided by uncertain chaotic sources~\cite{callegari05a}.

A discrete-time chaotic map, formed by the iteration of the output value in a transformation function, can be used for the generation of random numbers. Simple piecewise affine input-output (I/O) characteristics have been extensively used for the generation of random bits, e.g., the Bernoulli map \cite{callegari05a}, and the tent map~\cite{bean94}. The entropy source of a chaotic map is the inherent noise of the system, which is amplified in the positive gain feedback loop by the iteration of the output signal in the map function~\cite{callegari05a}. High Speed, capability of integration, and the high quality of the generated bits make the discrete-time chaotic maps very good candidates for high speed embeddable random number generation.

The practical application of both the tent map and the Bernoulli map can be hindered by noise and implementation errors, where they are unable to maintain the state of the system confined~\cite{callegari05a}. In this paper, we present a modified tent map that can be interchanged with the tent map in practical applications.

The rest of this paper is organized as follows. In Section~\ref{sec:map}, the fundamentals of discrete-time chaotic maps are reviewed, and the practical problems of the tent map are pointed out. In Section~\ref{sec:modified_tent}, we present the modified tent map and investigate its chaotic characteristics. In Section~\ref{sec:circuit}, we demonstrate the feasibility of implementing the presented modified tent map for true random number generation.
Section~\ref{sec:conclusion} concludes the paper.

\section{Tent Map: Fundamentals and History}
\label{sec:map}
In this section, after a brief introduction to discrete-time chaotic maps, we review the tent map and the previously proposed implementations. 
Discrete-Time Markov chaotic sources are a subclass of discrete-time chaotic nonlinear dynamical systems. A discrete-time chaotic system is formed by the iteration of the output signal through a transformation function $M(x)$ as given by
\begin{equation}
x_{n+1} = M(x_n) = M^n(x_0).
\end{equation}
In this equation, $n$ represents the time step, $x_0$ is the initial state of the system, and ${x_n}$ is the state of the system at time step $n$.

\begin{figure}
\vspace{0.2in}
\centering
\includegraphics[width=3.5in, angle=0]{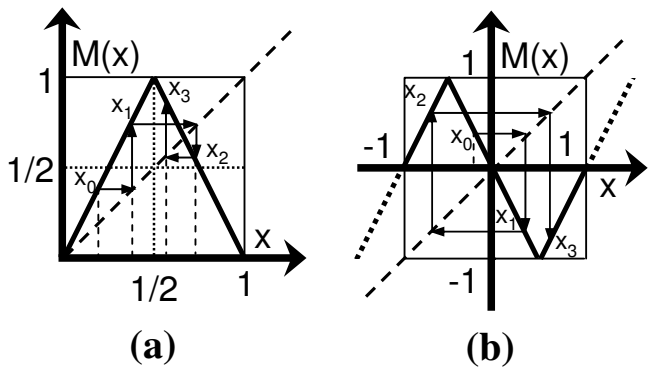}
\vspace{-0.08in}
\caption{(a)~The tent map. (b)~The presented modified tent map.}
\vspace{-0.15in}
\label{fig:maps}
\end{figure}

The tent map function $M(x): (0,1) \rightarrow (0,1)$, shown in Figure~\ref{fig:maps}~(a), is given by
\begin{equation}
M(x) =\left\{  \begin{array}{ll}
2x & 0<x\leq \frac{1}{2}, \\
2(1-x) &  \frac{1}{2}<x<1.
\end{array}
\right.
\label{eq:proposed_map}
\end{equation}
It can be shown that after several iterations, the density of states in the tent map asymptotically follows a uniform distribution, i.e., the state of the system is uniformly distributed in $(0,1)$ and the asymptotic density distribution $f(x)$ for the tent map satisfies $f(x)=1$~\cite{chaos_fractals}.

The Lyapunov exponent of a discrete-time map can be calculated from~\cite{chaos_dynamic}
\begin{equation}
\lambda = \int \ln  |M'(x)| f(x) dx,
\end{equation}
where $M'(x)$ is the derivative of the map function $M(x)$ and $f(x)$ is the asymptotic density distribution. A positive Lyapunov exponent implies the chaotic behavior in the system. The rate of the separation between the trajectories of very close initial states is given by $\lambda$. For tent map, $|M'(x)| = 2$ and $f(x)=1$, which result in $\lambda = \ln 2$.

Tent map circuit implementation has been proposed in~\cite{bean94, morie00}.
These circuit implementations suffer from the confinement problem in practice, i.e., the output value of the map can be trapped in a point outside the map due to noise or implementation errors~\cite{callegari97, callegari05a}. A tailed tent map has been presented in~\cite{callegari97} to solve the problem. The tailed tent map maintains the uniform asymptotic density distribution of the states, while not disturbing characteristics of the tent map. For example, the utilization of the tailed tent map for true random number generation results in the generation of correlated output binary sequence. 
In~\cite{callegari05b}, a hardware implementation has been proposed for the tent map based on reducing the slope, which will change the characteristics of the map and degrade the quality of the generated binary sequence in terms of the statistical characteristics.

\section{The Presented Modified Tent Map}
\label{sec:modified_tent}
In this section, we present the modified tent map and investigate its chaotic behavior. The presented modified tent map function $M(x): (-1,1) \rightarrow (-1,1)$, shown in Figure~\ref{fig:maps}~(b),
is given by
\begin{equation}
M(x) =\left\{  \begin{array}{ll}
2x+1 & -1<x\leq -\frac{1}{2},  \\
-2x & -\frac{1}{2}<x<\frac{1}{2}, \\
2x-1 &  \frac{1}{2}\leq x<1.
\end{array}
\right.
\label{eq:tent_map}
\end{equation}
In the modified tent map, the sign of the output value alternates in each iteration, i.e, $x_nx_{n+1}<0$.
Suppose $x_{1}^a$, $x_{2}^a$, and $x_{3}^a$ are the abscissa of three successive output values of the tent map, and $x_{1}^b$, $x_{2}^b$, and $x_{3}^b$ are the output values of the modified tent map. If $x_{1}^a = x_{1}^b$, we have $x_{2}^a= -x_{2}^b$ and $x_{3}^a= x_{3}^b$, as shown in Figure~\ref{fig:maps}~(a,b). In other words, for an equal initial state, the absolute value of the output sequence is equal for both the tent map and the modified tent map while the output of the modified tent map alternates between positive and negative values. Since the output values of the tent map and the modified tent map have equal absolute values, the Lyapunov exponent of the modified tent map is equal to that of the tent map, which was shown to be $\ln  {2}$.

\begin{figure}
\vspace{0.3in}
\centering
\includegraphics[width=3.1in, angle=0]{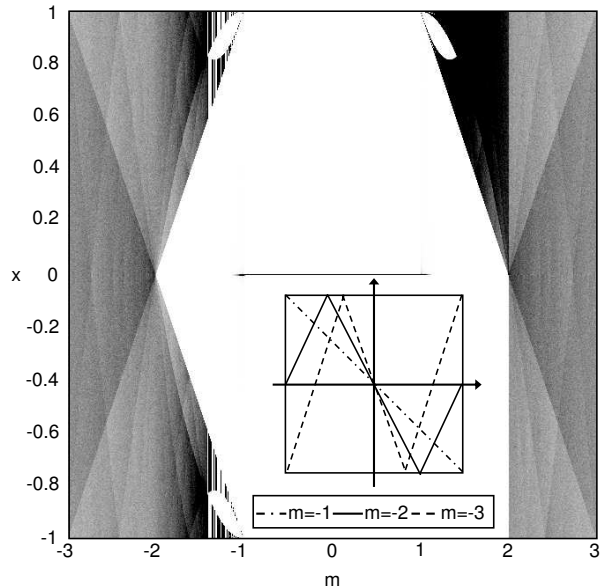}
\vspace{0.05in}
\caption{Bifurcation diagram of the modified tent map with bifurcation parameter m, defined in~(\ref{eq:bifurcation}). The generalized modified tent map with three different values of m ($-1$, $-2$, and $-3$) is shown in the inset.}
\vspace{-0.15in}
\label{fig:bifurcation}
\end{figure}

In order to further investigate the presented modified tent map, we generalize the map function by introducing a slope parameter $m$, where $m \in (-3,3)$, as given by
\begin{equation}
M(x) = \left\{
\begin{array}{ll}
-m(x +\frac{2}{|m|}) & -1<x \leq -\frac{1}{|m|}, \\
mx & -\frac{1}{|m|}<x<\frac{1}{|m|}, \\
-m(x- \frac{2}{|m|}) & \frac{1}{|m|} \leq x<1 .
\end{array} \right.
\label{eq:bifurcation}
\end{equation}
Here, $m=-2$ represents the modified tent map, as given by equation~(\ref{eq:tent_map}). Figure~\ref{fig:bifurcation} shows the bifurcation diagram of the generalized map function as a function of the bifurcation parameter $m$. A bifurcation diagram is obtained by following the output trajectory to find the possible long-term states of the system. Under certain conditions, the bifurcation diagram can also represent the density of the states due to the ergodicity~\cite{chaos_fractals, 758875}. In this diagram, the density of the states is proportional to the darkness, i.e. darker regions have higher density of states. This diagram is obtained assuming that the inherent initial noise is a very small positive value, i.e., $x_0 = 0^{+}$.

For $|m|>3$, the system could not maintain the state in the input range, which could not be used for chaotic applications. For $2<m<3$, the system is chaotic and the bifurcation diagram represents the steady state density distributions. Chaotic behavior can also be observed for $1<m\leq 2$. Since the initial state of the system is assumed to be positive, the output is confined in positive values.
For $|m|\leq 1$, the system does not demonstrate chaotic behavior. Regardless of the initial state of the system, the output will ultimately settle in zero. For $-2\leq m<-1$, the system also shows chaotic behavior. In this mode, the output alternates between positive and negative values with each iteration. In this region, the system does not have an asymptotic density distribution. For $-3<m<-2$, the system is chaotic and the bifurcation diagram represents its asymptotic density distribution.

For $m=2$, equation~(\ref{eq:bifurcation}) represents a map that is given by
\begin{equation}
M(x) =\left\{  \begin{array}{ll}
-2x-1 & -1<x\leq -\frac{1}{2},  \\
2x & -\frac{1}{2}<x<\frac{1}{2}, \\
-2x+1 &  \frac{1}{2}\leq x<1.
\end{array}
\right.
\label{eq:mirror_tent_map}
\end{equation}
As the equation shows, this is exactly like the tent map for $x>0$. Here, if $m$ is a little bit greater than $2$, the system could still maintain the state confined in the system since the output could jump to the negative side of the map, which is the fundamental difference between this map and the tent map. In the tent map, if the slope is more than $2$, the system could not maintain the state confined. Therefore, the tent map could not be used for practical applications under process variations and noise.

Although in the neighborhood of $m=-2$, the asymptotic density distribution does not exist, the modified tent map could be used in practical applications such as true random number generation based on the principles explained earlier.

\section{Implementation and Discussions}
\label{sec:circuit}
In this section, we propose the utilization of a current mode circuit implementation for the presented modified tent map. In current mode circuits, the signal is represented by the branch current instead of the node voltage in the voltage mode circuits. Current mode circuits have recently attracted great attention because of their capability of very low supply voltage operation~\cite{current_mode}. 
While scaling the feature size, the supply voltage is scaled more rapidly than the threshold voltage, which in turn decreases the overdrive voltage of the transistors and limits the voltage swing. In current mode circuits, high current gains can be achieved while the nodal voltages are floating~\cite{current_mode,roberts89}. 
Current mode circuits will be playing an important role in the future deep sub-micron technologies, e.g., telecommunications, analog signal processing, and multiprocessors~\cite{current_mode}.

\begin{figure}[t]
\vspace{0in}
\centering
\includegraphics[width=3.0in, angle=0]{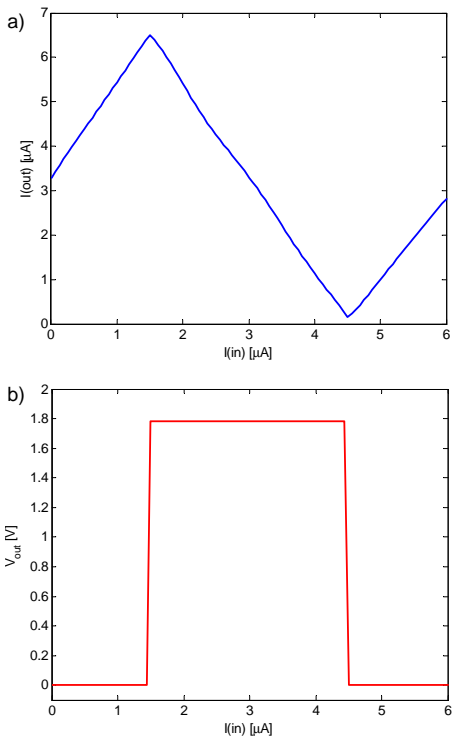}
\vspace{-0.08in}
\caption{(a)~The I/O current characteristics of the modified tent map. (b)~Output bit-generation pattern.}
\vspace{-0.15in}
\label{fig:map_bit_setti}
\end{figure}

The presented modified tent map is a \emph{continuous} piecewise affine map, where we implemented it by the affine interpolation of the value at the breakpoints, proposed in~\cite{758875}. In this method, the affine partitions are implemented using elementary blocks based on the detection of the breakpoints. Therefore, it is straight-forward to detect whether the output value is within a certain affine partition of the map, and generate a binary output sequence. We performed HSPICE simulations in TSMC $0.18~\mu m$ technology for the modified tent map circuit. The I/O characteristic curve of the map is shown in Figure~\ref{fig:map_bit_setti}(a), where the input and output range are equal to $6\mu A$. Investigating the transient response of the circuit with an input current pulse, we demonstrated that at the worst case the output would settle within $1\%$ of its final value in less than $25$~ns. Therefore, the circuit could be utilized at operation frequencies as high as $20$~MHz.

\begin{figure}[t]
\vspace{-1.3in}
\centering
\includegraphics[height=3.25in, angle=-90]{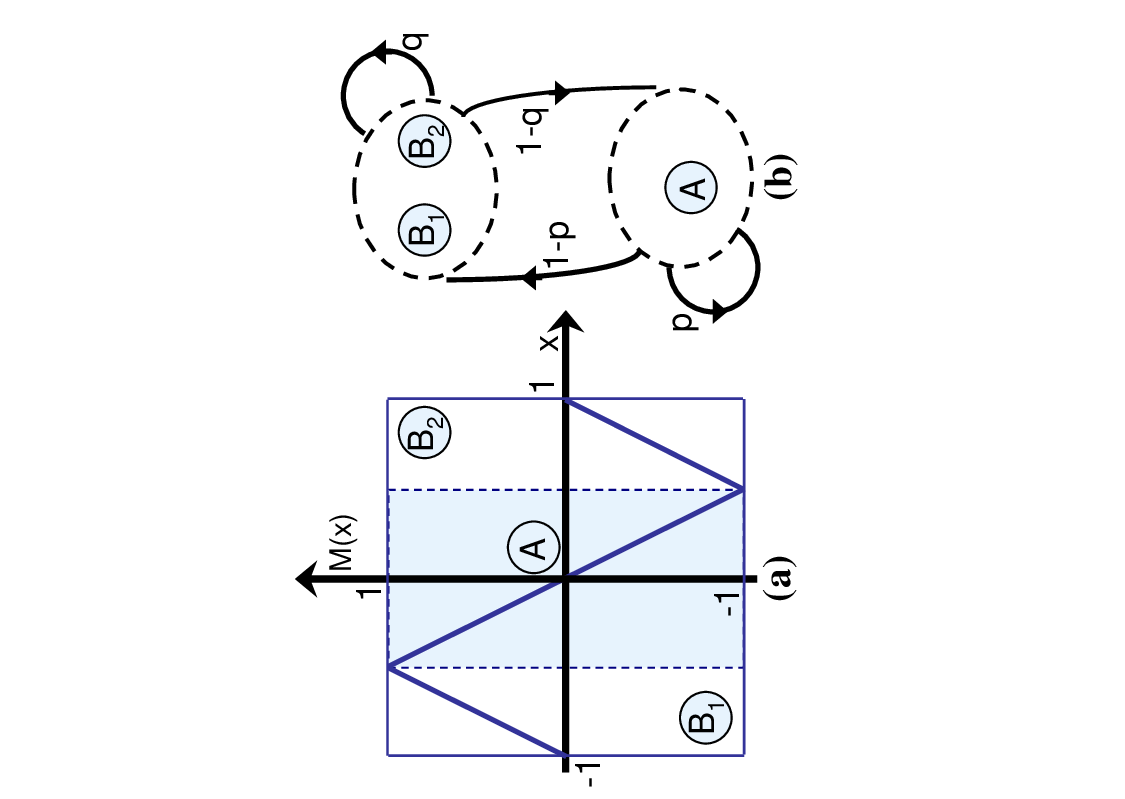}
\vspace{-1.3in}
\caption{(a)~The modified tent map and the corresponding affine partitions (b)~The Markov chain of states for true random number generation.}
\vspace{-0.2in}
\label{fig:map_ideal}
\end{figure}

We demonstrate the feasibility of binary generation from the presented modified tent map. In Figure~\ref{fig:map_ideal}(a), the partition $|x|<1/2$ is represented by state $A$, and the partitions $-1<x<-1/2$ and $1/2<x<1$ are represented by the states $B_1$ and $B_2$, respectively. In order to generate a binary output, the states $B_1$ and $B_2$ are combined in a macro-state $B$. The state $A$ corresponds to the bit $1$ and the state $B$ corresponds to the bit $0$. As discussed earlier, the binary sequence generated from a modified tent map is exactly the same as the binary sequence generated from a tent map due to the symmetry of the bit-generation process around the origin. Therefore, the bit-generation process can be modeled as a Markov chain similar to the tent map, as shown in Figure~\ref{fig:map_ideal}(b)~\cite{callegari05a}. In this Figure, $p$ is the probability of the generation of $1$ when the latest generated bit is $1$. $q$ is the probability of the generation of a $0$, when the previous bits is $0$. Since the functionality of the map is similar to the tent map, we have $p=q=1/2$ and the generated bits are truly random.  In Figure~\ref{fig:map_bit_setti}(b), the simulation result for the bit-generation circuit is presented, where it demonstrates very good agreement with the desired bit-generation pattern.

\section{Conclusion}
\label{sec:conclusion}
The tent map suffers from the inability to maintain the state confined in the input range of the map under noise and implementation errors. This fundamental issue limits the practical application of the tent map. In this paper, we presented a modified tent map, which can solve the state confinement problem of the tent map. The presented modified tent map could be used in all practical applications instead of the tent map. We demonstrated that the presented map could be used for true random number generation.

\balance
\end{document}